\documentstyle[11pt]{article}
\newtheorem{Theorem}{Theorem}[section]
\newtheorem{Lemma}[Theorem]{Lemma}

\def\df{\buildrel \rm def \over =} 
\newtheorem{Proposition}[Theorem]{Proposition}

\newtheorem{Example}[Theorem]{Example}

\def \rk{{\mbox rk}\,}
\def \dim{{\mbox dim}\,}
\def\df{\buildrel \rm def \over =} 
\def\PROOF{{\em Proof}\,: }              
\def\QED{~\hfill~ $\diamond$ \vspace{7mm}}
\def\ind{{\mbox ind}\,} 
\def\Ad{\,\mbox{Ad}\,}
\def\tr{{\mbox tr }}

\def\R{{\bf R}}
\def\C{\bf C}

\title{New examples of manifolds with completely integrable geodesic  
flows}

\author{G. P. Paternain\\ Department of Mathematics\\
SUNY Stony Brook\\Stony Brook, NY 11794 \and R. J. Spatzier
\thanks{Partially supported by the NSF and the Stony Brook Institute for
the Mathematical Sciences, Sloan Foundation Fellow and AMS Centennial
Fellow}\\ Department of Mathematics\\University of Michigan\\Ann Arbor, MI
48103}

\date{}

\begin{document}
\maketitle

{\em Dedicated to the memory of M. Strake}

\begin{abstract}
We construct Riemannian manifolds with completely integrable geodesic
flows, in particular various nonhomogeneous examples. The methods employed
are a modification of Thimm's method, Riemannian submersions and connected
sums.
\end{abstract}

\def \lie{g}
\def \ft{{
\beta}} 
\def \th{\bar{\theta}}

\section {\em Introduction}

A flow  $g_t$ on a symplectic $2n$-dimensional manifold $M$ is {\em 
completely integrable} if there are $n$ Poisson-commuting, $g_t$-invariant 
$C^{\infty}$-functions $f_1,\ldots,f_n$ whose differentials are 
independent a.e. in $M$.  Poincar\'{e} realized that complete 
integrability is an exceptional phenomenon.  Indeed, it was not 
until the past two decades that a large number of examples was 
discovered.  In this paper, we explore the narrower realm of 
geodesic flows.  Until now, very few examples of completely 
integrable geodesic flows were known.  The classical examples are 
the flat tori, surfaces of revolution (Clairaut), $n$-dimensional 
ellipsoids with different principal axes (Jacobi) and $SO(3)$ with a 
left invariant metric (Euler). More recent examples are semisimple 
Lie groups with certain left invariant metrics due to Mishchenko and 
Fomenko \cite{Mish}. Then Thimm devised a new method for constructing
first integrals in involution on homogeneous spaces \cite{Thimm}.   In
particular, he proved  the complete integrability of  the geodesic flow
on real    
and complex Grassmannians. Guillemin and Sternberg 
conceptualized this method and found further examples \cite{G2}. The 
spaces  obtainable in this way have essentially been classified by 
Kr\"{a}mer \cite{Kr}.  For results concerning  the complete integrability
of the geodesic flows of other symmetric spaces of the compact type  we
refer to \cite{Fomenko}.

In this paper, we exhibit several new examples of Riemannian 
manifolds with completely integrable geodesic flows, and in 
particular various nonhomogeneous examples. We use several new 
techniques in these constructions.

     The first construction is a simple variation of Thimm's method. 
In his method, the moment map of a Lie group action is used to pull 
back a family of Poisson-commuting functions from the Lie algebra to 
the symplectic space in question. While Thimm considered the case of 
an action of a group by isometries on a homogeneous space, one can 
instead use the isometry group together with the geodesic flow. This 
generalizes the construction of integrals on surfaces of revolution. 
As simple as it is, this variation of Thimm's method already yields 
several new examples such as the Wallach manifold $SU(3)/T^2$ where 
$T^2$ is a maximal torus in $SU(3)$. We explore this in Section 2. 
 
In Section 3, we study the symplectic structure of a Riemannian 
submersion. When the submersion is given by an action of a Lie group 
$G$, the tangent bundle of the base space is symplectomorphic with 
the Marsden-Weinstein reduction of the tangent bundle of the total 
space with respect to the action of the group by derivatives. In 
particular, we see that $G$-invariant Poisson-commuting functions 
descend to Poisson-commuting functions. Thus the base space of the 
submersion has completely integrable geodesic flow if enough $G$-
invariant Poisson-commuting functions descend to independent 
functions. This is the essence of the {\em submersion method}. 

     In Section 4 we construct various examples using the submersion 
method. To apply it, we observe that sometimes the integrals arising 
>from the Thimm method are invariant under the action of a subgroup 
of the isometry group. Then we show in the examples that enough 
Poisson-commuting functions descend to independent functions on the 
quotient space. Unfortunately, the independence of the functions on 
the base space is far from automatic.  However we     show that in general
if $X$ has completely integrable geodesic flow  
and admits an $S^1$-action that leaves the integrals invariant and  
$N$ is a surface of revolution then the geodesic flow of $X \times 
_{S^1} N$ is completely integrable.  Particular examples are ${\bf 
CP} ^n \#- {\bf CP} ^n$ and surface bundles over the  Eschenburg examples.
In the latter the base space is a quotient space  
of $SU(3)$  by  $S^1$-action acting both from the left and the right. 
Some of these surface bundles are known to be strongly inhomogeneous, that
is, they  
do not have the homotopy type of a compact homogeneous space 
\cite{Spatzier-Strake}.  Next, we show that certain Eschenburg examples 
themselves have completely integrable geodesic flows. Again most of 
these spaces are strongly inhomogeneous. They do not fall under the 
general submersion example above. Rather, we use the submersion 
method directly, and establish independence of sufficiently many 
functions by explicit computation. Note also that while the isometry
 groups of the Eschenburg manifolds are non-trivial, they are not big
enough for the Thimm method to apply, due to dimensional reasons. Finally,
we show that the geodesic flow of the  exotic sphere used by Gromoll and
Meyer in \cite{Gromoll-Meyer} is completely integrable. Here the integrals
come both from a Thimm construction combined with the submersion method as
well as from the isometry group.  Let us remark that the geodesic flows of
certain Kervaire spheres also admit a complete set of integrals on an open
dense subset of the tangent bundle. It is not clear however whether these
integrals extend to the full tangent bundle.         

In Section 5, we use a glueing technique to construct metrics with 
completely integrable geodesic flows on  ${\bf CP} ^n \# {\bf CP} 
^n$ for $n$ odd. 

The second author is grateful to M. Strake who had introduced him to the
Eschenburg examples. It also was first in discussions with him that the
possibility of the complete integrability of the geodesic flows of the
Eschenburg examples arose. We  also thank D. Gromoll and B. Kasper for
helpful comments and discussions.

\section{\em Variations on the Thimm method.}

First we recall Thimm's  construction as modified by Guillemin and  
Sternberg \cite{G1,G2,Thimm}.  
We refer to \cite{G2} for more details. 

Let $N$ be a symplectic  space with a Hamiltonian action of a Lie  
group $G$. Such an   action is called  {\em multiplicity free} if  
the algebra of the $G$-invariant functions on $N$ is commutative  
under the Poisson bracket \cite[p. 361]{G1}. Let  
$\Phi : N \rightarrow g^{*} $ denote the moment map of the action. Let 
 $\{1\}=G_{l}\subset G_{l-1}\subset ...\subset G_{1}=G$  
be an ascending chain of Lie subgroups of $G$, and denote their Lie  
algebras by $g _i$. Furnish each coadjoint orbit of $G_i$ in  
$g _i ^*$ with the Kostant-Kirillov symplectic structure. Then  
each subgroup $G _{i+1}$ acts on each orbit in $g_i ^*$ in a  
Hamiltonian way.  The moment maps are just the restrictions of the  
dual maps $j_i : g _i ^* \rightarrow g_{i+1} ^*$ to the  
coadjoint orbits. We will call the chain $G_i$ {\em multiplicity-
free} if the actions of the $G _{i+1}$
on the coadjoint orbits of $G_i$ on $g_i ^*$ are multiplicity  
free. This is quite a strong condition on the chain $G_i$. 
For compact groups it forces the $G_i$ to be locally isomorphic to
$SO(n)$, $SU(n)$, tori or products of these \cite{Heckman,Kr2}.

If the $G_i$ are a multiplicity-free chain and the action of $G$ on 
$N$ is multiplicity-free, then any $G$-invariant Hamiltionian on $N$ 
is completely integrable \cite[p. 366]{G1}. This is the
essence of the {\em Thimm method}. This setup was studied in detail in 
\cite{Thimm,G2,G3}. If $N$ is the cotangent bundle of a manifold $M$ and
$G$ acts by derivatives then $M$
is a homogeneous space $G/K$ \cite{G3}. In this case, one calls 
$(G,K)$ a {\em Gelfand pair}. They have been classified by Kramer in 
\cite{Kr}.

We observe that a variation of the Thimm method also gives
complete integrability of some geodesic flows on homogeneous spaces $G/K$
even when the pair $(G,K)$ is not a Gelfand pair.
First let \ind $G$ denote the index of $G$. It is 
defined as the codimension of a generic orbit of the coadjoint action of
$G$ on $l^{*}$ (if $G$ is semisimple \ind $G$=\rk $G$).
Consider now a homogeneous space $G/K$ that verifies the following
conditions:
\[\begin{array}{l}
(1) \mbox{ dim }  G=2 \,\dim K \,+\, ind G\,+\,2 \\
(2)  \mbox{  the isotropy group of  }  G  \mbox{ at  some } v \in T(G/K)
 \mbox{    has dimension zero.} 
\end{array}\]

%\vspace{.20in}

Denote by $R$ the Hamiltonian action generated by the geodesic flow of
some left invariant metric. Let $\hat{G}=G\times R$. Clearly $\hat{G}$
acts
by Hamiltonian transformations and leaves the quadratic form associated
with the metric invariant.
Suppose the left invariant metric on $G/K$ has a geodesic which is not the
orbit of a 1-parameter subgroup of $G$.
Then a.e. the isotropy of $\hat{G}$ has dimension zero and
\dim$\hat{G}$+\ind$\hat{G}$=\dim$T(G/K)$.
By a dimension count, we deduce that the isotropy groups of the coadjoint
action of $\hat{G}$ act transitively on the regular level surfaces of the
moment map of $\hat{G}$.  
This implies that the action of $\hat{G}$ is multiplicity free by the
equivalences stated in \cite{G2}.
Therefore the Thimm method can be applied whenever 
a multiplicity free chain can be constructed for $G$.

Now, let us rewrite
condition (2).
Assume $G$ is compact and denote by $(\;,\;)$ some bi-invariant metric.
Let $k$ denote the Lie algebra of $K$ and $k^{\bot}$ the orthogonal
complement  with respect to $(\;,\;)$. 
Then it is easy to check that (2) is equivalent to:
\[\begin{array}{l}
(3) \mbox{ For some } X\in k^{\bot} \mbox{ we have } \dim [X,k] =\dim
k.\hfill
\end{array}\]
%\vspace{.20in}

Let us see some examples:
\begin{Example}{\rm
 Consider the homogeneous space $SU(3)/T^{2}$ where $T^{2}$ is a maximal
torus in $SU(3)$. This manifold can be also considered as the space of
flags in ${\bf CP^{2}}$. Since \dim$SU(3)=8$ and \ind$SU(3)=2$ condition
(1) is clearly verified.
 We will check condition (3).

The Lie algebra of $SU(3)$ consists of all the skew hermitian matrices
with trace zero. In this case $k$ consists of all the matrices $Y$ of the
form:
\[Y=\left(\begin{array}{clc}\alpha&0&0\\0&\beta&0\\0&0&
	\gamma\end{array}\right)\]
where $\alpha$, $\beta$ and $\gamma$ are purely imaginary and their sum is
zero.

Consider the Killing metric on $SU(3)$ i.e.
$(X,Y)=-\frac{1}{2}Re\;tr(XY)$.
 With respect to this product $k^{\bot}$ is the subset of $su(3)$ given by
the matrices with zero entries on the diagonal. Let $X\in k^{\bot}$ be
given by :
\[X=\left(\begin{array}{ccc}0&1&1\\-1&0&0\\-1&0&0\end{array}\right)\]
Take $Y\in k$ as before and compute $[X,Y]$. We get:
\[[X,Y]=\left(\begin{array}{ccc}0&\beta-\alpha&\gamma-\alpha\\
\beta-\alpha&0&0\\\gamma-\alpha&0&0\end{array}\right)\]
Then we clearly have \dim$[X,k]$=\dim$k$=2 and condition (3) is verified.

Therefore the geodesic flow of a left invariant metric on $SU(3)/T^{2}$ is
completely integrable provided that not every geodesic is the orbit of a
1-parameter subgroup of $SU(3)$.

Exactly the same arguments can be applied to other spaces.   In
particular, the geodesic flow of a left invariant metric on
$SO(n+1)/SO(n-1)$ is completely integrable provided that not every
geodesic is the orbit of a 1-parameter subgroup of $SO(n+1)$. The complete
integrability of the geodesic flow of the normal homogeneous metrics on
$SO(n+1)/SO(n-1)$ was obtained by Thimm \cite[Proposition 5.3]{Thimm}.
Here the original Thimm method works since the natural action of
$SO(n+1)\times SO(2)$ on the tangent bundle of $SO(n+1)/SO(n-1)$ is
multiplicity free.}

\end{Example}

%% GABRIEL, PLEASE INSERT YOUR  OBSERVATION ABOUT THE CODIM 1 ACTIONS

\section{\em Submersion metrics and reduced spaces}

An especially nice class of Riemannian submersions is that given by
isometric group actions. Their main symplectic feature, as we will see, is
that their tangent bundles are  Marsden-Weinstein reductions of the
tangent bundles of the  total spaces. This and other basic symplectic
properties are fundamental to the examples studied in the  remaining
sections. We refer to \cite[ch. 9]{Besse} and \cite[section 26]{G1} for
all the basic  definitions.

Let a Lie group $G$  of dimension $m$ act on a Riemannian manifold $M$
with metric  
$\langle\, , \, \rangle _M$ by isometries. We endow the tangent bundle
$TM$ of $M$  with the symplectic structure $\bar{\omega}$ obtained by
pulling back the canonical  symplectic structure on the cotangent bundle
$T^* M$ by the metric. Then $G$ acts  symplectically on $TM$. Let $\lie$
denote the Lie algebra of $G$. 

\begin{Lemma}               
The moment map $\Phi: TM \rightarrow \lie ^*$ is given by
\[ \Phi(v) (\xi) = \langle v, \xi (\ft (v)) \rangle \]
for $v \in TM$ and $\xi \in \lie$ where $\ft$ maps a tangent vector to  
its foot point. 
\end{Lemma}

\PROOF Recall the formula for the moment map on the cotangent  
bundle, namely $\Phi^* (v) = v (\xi (\ft (v)))$ for $v \in T^* M$ and  
$\xi \in \lie$ \cite[p. 222]{G1}. This readily  
implies the claim since the symplectic structure on $TM$ is the  
pullback under  the Riemannian structure.
\QED

Now suppose that $G$ acts on $M$ without isotropy. Set $B = M/G$ and  
endow $B$ with the submersion metric. Denote  by $\pi: M \rightarrow  
B$ the quotient map. 

\begin{Lemma}
The moment map intersects the trivial  coadjoint orbit $\{0\}$ in $\lie
^*$  cleanly, i.e. $\Phi^{-1}(0)$ is a submanifold of $TM$ and at each
point $x \in \Phi^{-1}(0) $ we have $T _x (\Phi^{-1}(0))  = d\Phi ^{-1} (T
_0 (0)) = d\Phi ^{-1} (0)$. Moreover, $\Phi^{-1} (0) $ is the set of all
horizontal vectors.
\end{Lemma}

\PROOF First note that $\Phi^{-1} (0)= \{v \in TM \mid \mbox{for all }\xi
\in \lie \: \langle v, \xi(\ft ( v))\rangle =0 \} $ is the  
set of all horizontal vectors, and thus a manifold. 

For  $w = (w_1,w_2) \in T _x (\Phi^{-1}(0))$ let $\{p_t\} \subset M$ and
$\{ v _t\} \subset TM$ be $C^1$-paths such that $w_1 =\frac{d}{dt} \Big|
_{t=0} p_t $ and $w_2 =\frac{d}{dt} \Big| _{t=0} v_t  $.  
Let $ \tilde{\xi _1}, \ldots,\tilde{ \xi _m}$ be a basis
 for $\lie$. Choose a coordinate system for $M$ about $x$ such that the
first $m$ coordinates are given by $\xi _1 = \tilde{\xi _1} (y), \ldots,
\xi _m = \tilde{\xi _m }(y)$. In this coordinate system we may write $v_t
= h_t + \eta + t \xi = h_t + \tilde{\eta}(p_t) + t \tilde{\xi} (p_t)$
where the $h_t$ are horizontal
 and $\tilde{\eta}$ and $\tilde{\xi} \in \lie$. Since $\lie ^*$ is a
vector space, $T^* \lie$ is canonically identified with $\lie ^*$ and we
have 
\[  (\frac{d}{dt} \Big| _{t=0} c_t )(\zeta) \df \frac{d}{dt} \Big| _{t=0}
(c_t (\zeta))\] 
for $\zeta \in \lie ^*$ and $c_t$ a $C^1$-path in $\lie ^*$.  
Since the $h_t$ are horizontal we get  
\[ d\Phi _x (w) (\xi) = \frac{d}{dt} \Big| _{t=0} \langle v_t, \xi (p_t)
\rangle = \frac{d}{dt} \Big| _{t=0} \langle  h_t + \eta + t \xi (p_t), \xi
(p_t) \rangle = \langle \xi (p_t), \xi (p_t) \rangle .\]
Thus $d\Phi _x (w) =0$ if and only if $ \xi =0$ or equivalently if $w \in
T _x (\Phi^{-1}(0))$.
\QED

By the last lemma the Marsden-Weinstein reduced space $TM / / G$ of $TM$
with respect to the $\{0\}$-coadjoint orbit is defined. Recall that it can
be identified with the reduced space $\Phi^{-1} (0) / G$ \cite[p.
192]{G1}.

\begin{Proposition} 
The Marsden-Weinstein reduced space $TM / / G$ with respect to the
$\{0\}$-coadjoint orbit  is symplectomorphic with $TB$.
\end{Proposition}

 \PROOF  Since  $\Phi^{-1} (0)$ is the  
set of all horizontal vectors, $TM / / G$ is diffeomorphic  
with $TB$. We will identify the two henceforth.

The canonical symplectic structure on $T^* B$ pulled back to $TB$ by   the
Riemannian metric $\langle\, , \, \rangle _B $ is the 2-form   $\omega =
-d \theta$ where $\theta $ is the 1-form defined by  
$\theta _x (v) = \langle x, d(\ft) _x v \rangle _B$ for $x \in TB$  
and $v \in T_x TB$. Similarly denote by $\th$  the  1-form on  $TM$, with
$\bar{\omega} = -d \th $. Let  
$\bar{x} \in \Phi^{-1} (0)$ and $\bar{v} \in T_{\bar{x}} \Phi^{-1} (0)$.  
Then $\bar{x}$ and $d (\ft) _{\bar{x}} \bar{v}$ are horizontal  
vectors, and therefore
\[\th _{\bar{x}} (\bar{v}) = \langle \bar{x}, d(\ft) _{\bar{x}}  
\bar{v} \rangle _M =   \langle d \pi( \bar{x}),d \pi( d(\ft)  
_{\pi \bar{x}} \bar{v}) \rangle _B = \theta _{\pi \bar{x}} (d \pi  
\bar{v}).\]
If we restrict $d \pi$ to $\Phi^{-1} (0)$, we see  that $d \pi ^* \theta =
\th \mid _{\Phi^{-1} (0)}$. Hence $\bar{\omega} \mid _{\Phi^{-1} (0)} =  
d \pi ^* \omega$. By definition, the symplectic form $\omega _r$ on   the
reduced space  satisfies $d \pi ^* \omega _r = \bar{\omega}  
\mid _{\Phi^{-1} (0)}$. Since $d \pi$ is surjective, we get $\omega =  
\omega _r$.
\QED

For a $C^1$-function  $f$ on a symplectic manifold $X$ we denote by
$\xi_f$ the associated {\em Hamiltonian vector field}. Suppose $G$
acts on $X$ in a Hamiltonian way with moment map $\Phi: X \rightarrow
\lie ^*$  that intersects $\{0\}$ in $\lie ^*$ cleanly.
If  $f$ is  $G$-invariant then $\xi _f$ is tangent to $\Phi^{-1} (0)$.
Let $\rho $ be the projection from $\Phi^{-1} (0)$ to the reduced
space $Y = \Phi^{-1} (0) /G$. If $f':Y \rightarrow \R$ denotes the
induced function then $\xi _{f'} = d \rho (\xi _f)$ \cite[Appendix
5C]{Arnold}.
 
\begin{Lemma}
Let $f,g: X \rightarrow \R$ be $G$-invariant functions on $X$ and
$f',g':Y \rightarrow \R$ the induced functions. Then we have 
$$\{f',g'\} \circ
\rho = \{f,g\} \mid _{\Phi^{-1} (0)}.$$
 In particular, Poisson-commuting
 $G$-invariant functions descend to Poisson-commuting functions.
\end{Lemma}
 
\PROOF This follows from the discussion before the lemma and
\[\xi_{ \{f',g'\}}  = [\xi _{f'},\xi _{g'}]  = d \rho  \,([ \xi _f, \xi
_g].\]
\QED

\section{\em Submersion examples.}
The idea of the constructions below is that sometimes the integrals that
arise from the Thimm method are invariant under a subgroup of the isometry
group.
 Then one can construct integrals for the quotient space by this subgroup
endowed with the submersion metric. Of course, the main problem is to show
independence of the integrals thus obtained.

We need to describe the Thimm method in more detail to understand the
invariance properties of the Thimm integrals. First one finds a maximal
family of functions on $g^{*}$ in involution which are functionally
independent.
 Their construction is inductive: one pulls back a family of such
functions already constructed on $g^{*}_{i+1}$ by $j_{i}$ and appends a
maximal number of functionally independent $G_{i}$-invariant functions on
$g_{i}$.
 Now one can pull back this family of functions on $g^{*}$ to the
symplectic space $N$ using the moment map $\Phi$.
 Under our hypothesis, we get $n:=\dim N/2$ many functions
$f_{1},...,f_{n}$ in involution on $N$ which almost everywhere are
functionally independent.
 Furthermore they commute with any $G$-invariant Hamiltonian on $N$.

Next we observe some invariance properties of these integrals.

\begin{Lemma}
Let $j:g'\rightarrow g$ be a Lie subalgebra corresponding to the subgroup
$G'$ of $G$.
Let $\phi':g'^{*}\rightarrow\R$ be a function invariant under the
coadjoint action of $G'$. Let $\phi=j^{*}(\phi ')=\phi ' \circ j^{*}$.

{\bf a)} If $\tau$ is in the centralizer of $g'$ in $G$ then $\phi$ is
invariant under $\tau$.

{\bf b)} If $\tau\in G'$ then $\phi$ is invariant under $\tau$.
\end{Lemma}

\PROOF
Suppose $\tau$ is in the normalizer of $g'$ in $G$.
 Then 
\[\tau(\phi)=\phi\circ \Ad^{*}(\tau)=\phi '\circ j^{*}\circ
\Ad^{*}(\tau)=\phi '\circ \Ad^{*}(\tau)\circ j^{*}.\]
If $\tau$ centralizes $g'$ then $\Ad^{*}(\tau)=id\mid_{g'}$, and {\bf a)}
follows. 
 For {\bf b)} recall that $\phi '$ is $G'$-invariant.
\QED

In particular we see that any function $f_{i}$ constructed above is
$G_{l-1}$-invariant.
 In fact, since the moment map is equivariant, it suffices to see this for
the functions $\phi_{i}: l^{*}\rightarrow \R$ constructed above.
 Since $G_{l-1}$ is contained in all the other subalgebras (except
$\{1\}$) this follows from Lemma 4.1 b).
\vspace{1.5em}

In the remainder of this section, we combine these observations with the
results of the previous sections.
 We first describe a general construction of new manifolds supporting
integrable geodesic flows from known ones.

Let $X$ be a complete Riemannian manifold of dimension $n$ whose geodesic
flow is completely integrable.
 Suppose $X$ admits a free $S^{1}$ action by isometries that leaves the
integrals invariant.
 Let $N$ be a complete surface of revolution and consider the diagonal
action by isometries on $X\times N$. Since this action is free we can
consider the quotient manifold $M=X\times_{S^{1}}N$ and endow it with the
submersion metric.
  Then we have:

\begin{Proposition}
 The geodesic flow on $M$ is completely integrable.
\end{Proposition}

\PROOF
Denote by $f_{1},...,f_{n}$ the integrals coming from $X$ and by $g$ the
metric on $N$.
 These functions extend to $T(X\times N)=TX\times TN$ in the obvious way
and their extensions will be denoted by
$\hat{f_{1}},...,\hat{f_{n}},\hat{g}$. 
They clearly are integrals of the geodesic flow given by the product
metric on $X\times N$.
 Since all these integrals are invariant under $S^{1}$ by hypothesis, they
descend to $\Phi^{-1}(0)/S^{1}=TM$ where $\Phi$ is the moment map
corresponding to the diagonal action of $S^{1}$ on $X\times N$.
 Hence we get $n+1$ integrals in involution for the geodesic flow on $M$
(cf. Proposition 3.3 and Lemma 3.4). 
 We only need to show that they are functionally independent almost
everywhere.

For this consider a point $(p_{1},p_{2})\in X\times N$. Set
$H=\Phi^{-1}(0)\cap T_{(p_{1},p_{2})}(X\times N)$ which is nothing but the
set of horizontal vectors at $(p_{1},p_{2})$.
 Let $\tau :H\rightarrow T_{p_{1}}X$ be the restriction of the projection
map.
 Clearly $\tau$ is onto if $S^{1}$ does not fix $p_{2}$.
 For a.e. $p_{1}$ and a.e. $v_{1}\in T_{p_{1}}X$ the vector fields
$\xi_{f_{1}},...,\xi_{f_{n}}$ at $v_{1}$ are linearly independent by
assumption. Take such a $v_{1}$ and let $v_{2}$ be a non-zero vector on
the projection of $\tau^{-1}(v_{1})$ over $T_{p_{2}}N$.
 Moving $v_{1}$ a little if necessary, we can choose $v_{2}$ so that the
geodesic through $v_{2}$ is not an orbit of the $S^{1}$ action.
 We will now check independence at the projection of $(v_{1},v_{2})$ to
$\Phi^{-1}(0)/S^{1}$.

It is clear that $\xi_{\hat{f_{1}}},...,\xi_{\hat{f_{n}}},\xi_{\hat{g}}$
are linearly independent at $(v_{1},v_{2})$.
 Therefore we need to show that the tangent vector field $W$ to the orbits
of $S^{1}$ on $\Phi^{-1}(0)$ at $(v_{1},v_{2})$ does not belong to the
vector space spanned by
$\xi_{\hat{f_{1}}},...,\xi_{\hat{f_{n}}},\xi_{\hat{g}}$.
 Let $W_{1}$ denote the tangent field to the orbits of $S^{1}$ on $TX$ and
$W_{2}$ the corresponding tangent field to the orbits of $S^{1}$ on $TN$.
 Then $W$ can be written as $(-W_{1},W_{2})$ where $W_{1}(v_{1})\in
T_{v_{1}}TX$ and $W_{2}(v_{2})\in T_{v_{2}}TN$.
 Observe now that $\xi_{\hat{f_{i}}}\in TTX$ and $\xi_{\hat{g}}\in TTN$.
 Hence if $W$ belongs to the space spanned by the $\xi_{\hat{f_{i}}}\;'s$
and $\xi_{\hat{g}}$ at $(v_{1},v_{2})$ we deduce that $W_{2}$ and
$\xi_{\hat{g}}$ are collinear at $v_{2}$. This implies that the geodesic
through $v_{2}$ is an orbit of the $S^{1}$-action.
 This is a contradiction to the choice of $v_{2}$.
\QED

Let us see some applications of the previous proposition.

\begin{Example}{\rm
Consider $SO(n)$ endowed with its standard bi-invariant metric.
 Then $SO(n)$ can be viewed as a symmetric space of $SO(n)\times SO(n)$ in
the usual way.
 Consider the ascending chain of subgoups:
\[\{1\}\subset SO(2)\times\{1\}\subset ...\subset SO(n)\times\{1\}\subset
...\subset SO(n)\times SO(n)\]
This chain as well as the action of $SO(n)\times SO(n)$ on $TSO(n)$ are
multiplicity free and  the hypotheses of the Thimm method hold. Thus we
recover the well-known fact that the geodesic flow on $SO(n)$ is
completely integrable. By Lemma 4.1 the integrals are all invariant under
$SO(2)\times\{1\}=S^{1}$.
Now the proposition implies that $M=SO(n)\times_{S^{1}}N$ endowed with the
submersion metric supports a completely integrable geodesic flow.
 Similar arguments apply to the case of $SU(n)$.

Let us describe the manifolds we get for the special case of
$SU(2)=Spin(3)=S^{3}$ and $M^{4}=SU(2)\times_{S^{1}}N^{2}$.
  If $N$ is the 2-sphere then $M$ is the non-trivial $S^{2}$-bundle over
$S^{2}$ which is diffeomorphic to ${\bf CP^{2}\#-CP^{2}}$. It is known
that $M$ is not diffeomorphic to any homogeneous space (cf. \cite{C}). 
 If $N$ is euclidean 2-space, then $M$ is the normal bundle of $\C P^{1}$
in $\C P^{2}$. Moreover if we consider $S^{1}$ acting on the plane by
rotating n times then the corresponding $M_{n}$ give all the line bundles
over $\C P^{1}$.
 Finally we want to point out that the metric on $SU(2)$ does not need to
be bi-invariant. For any left invariant metric the arguments work because
the action of $\hat{G}=SU(2)\times R$ as in Section 2 is multiplicity
free.}
\end{Example}

\begin{Example}{\rm Next we will construct a large class of metrics on $\C
P^{n}\#-\C P^{n}$ with completely integrable geodesic flows, generalizing
the last example.

Consider the Hopf fibration $S^{1}\rightarrow S^{2n+1}\rightarrow \C
P^{n}$.
 Denote by $g_{t}$ the metric on $S^{2n+1}$ which is obtained from the
standard metric by multiplying with $t^{2}$ in the directions tangent to
the $S^{1}$-orbits.
The canonical action of the group $SU(n+1)$ on $S^{2n+1}$ is by isometries
and  commutes with the $S^{1}$ action. Hence the group $SU(n+1)\times
S^{1}$ acts on $S^{2n+1}$ by isometries.
 It is known that $(S^{2n+1},g_{t})$ can be viewed as distance spheres on
$\C P^{n+1}$ with the metric induced by the Fubini-Study metric.
 For $t\leq \frac{n+1}{2n}$ they are called Berger spheres. We refer to
\cite{Z} for details.
 The action of $SU(n+1)\times S^{1}$ on $TS^{2n+1}$ is multiplicity free.
Choosing a suitable chain it follows from the Thimm method that the
geodesic flow on $(S^{2n+1},g_{t})$ is completely integrable. By Lemma 4.1
the integrals are all invariant under the $S^{1}$-action.
 Then Proposition 4.1 shows that the geodesic flow on
$M=S^{2n+1}\times_{S^{1}}N$ is completely integrable for all real $t$.
 If we take $N=S^{2}$ then the corresponding $M$ is diffeomorphic to $\C
P^{n+1}\#-\C P^{n+1}$. 
 If $N$ euclidean 2-space then $M$ is the normal bundle of $\C P^{n}$ in
$\C P^{n+1}$.}
\end{Example}

\begin{Example}{\rm Let $G_{n-1,2}({\bf R})=SO(n+1)/SO(n-1)\times SO(2)$
denote the Grassmannian of 2-planes in n+1-space.
 Consider the fibration $S^{1}\rightarrow SO(n+1)/SO(n-1)\rightarrow
G_{n-1,2}({\bf R})$, where $S^{1}$ acts on $SO(n+1)/SO(n-1)$ by right
translations.
 As  mentioned at the end of Section 2, the action of $SO(n+1)\times
S^{1}$ on the tangent bundle of $SO(n+1)/SO(n-1)$ is multiplicity free.
Consider metrics $g_{t}$ on $SO(n+1)/SO(n-1)$ obtained from the normal
homogeneous metric by multiplying with $t^{2}$ in the directions tangent
to the $S^{1}$-orbits.
 Thus we can argue as in Example 4.4 to deduce that the geodesic flow on
$M=SO(n+1)/SO(n-1)\times_{S^{1}}N$ is completely integrable for all real
$t$.
 If $N$ is the 2-sphere, $M$ is a sphere bundle over the Grassmannian
$G_{n-1,2}({\bf R})$.}
\end{Example}

\begin{Example}{\rm Next consider surface bundles over the so called
Eschenburg examples \cite{Eschenburg} (we will discuss the Eschenburg
examples themselves below).
 Consider the group $SU(3)$ with its standard bi-invariant metric and let
$SU(3)\times SU(3)$ act on $SU(3)$ by $(g_{1},g_{2})x=g_{1}xg_{2}^{-1}$.
 Let $k,l,p,q$ be a set of relatively prime integers.  Define a
one-parameter subgroup of $SU(3)\times SU(3)$ by
\[U_{klpq}=\{\exp 2\pi
it(\mbox{diag}\,(k,l,-k-l),\mbox{diag}\,(p,q,-p-q))\mid t\in \R\}.\]
For certain choices of $k,l,p$ and $q$ the action of $U_{klpq}$ on $SU(3)$
is fixed point free, in particular for the quadruple $(1,-1,2m,2m)$
\cite[Proposition 21]{Eschenburg}. 

Consider the ascending chain of subgroups:
$\{1\}\times U(1)\subset U(1)\times U(1)\subset U(1)\times  U(2)\subset
U(2)\times U(2)\subset U(2)\times SU(3)\subset SU(3)\times SU(3)$ 
where $U(1)$ and $U(2)$ are embedded into $SU(3)$ by adjusting the
$(3,3)$-entry in the matrix in the obvious way. 
 Note that $(1,\exp 2\pi it\,\mbox{diag}\, (2m,2m,-4m))$ and $(\exp 2\pi
it\, \mbox{diag} \,(1,-1,0),1)$ either belong to or centralize any
subgroup in this chain.
 Thus all the first integrals on $TSU(3)$ are invariant under these
one-parameter subgroup and thus under $U_{1,-1,2m,2m}$.

>From the last proposition  we deduce that
$M_{m}=SU(3)\times_{U_{1,-1,2m,2m}}N$ endowed with the submersion metric
supports a completely integrable geodesic flow .
 If $N$ is the 2-sphere, $M_{m}$ is a sphere bundle over Eschenburg's
strongly inhomogeneous 7-manifold $SU(3)/U_{1,-1,2m,2m}$. These spaces
where studied in \cite{Spatzier-Strake}. Metrically they have higher rank
and topologically are strongly inhomogeneous  and irreducible
(\cite[Proposition 4.2 and 4.6]{Spatzier-Strake}).}
\end{Example}

Finally let us study some submersions that do not have the product type
used in Proposition 4.2. The observations concerning the invariance of the
Thimm integrals from the beginning of this section however are still
crucial.  Unfortunately, the calculations necessary become much more
complicated. 

\begin{Example}{\rm  Here we will study the Eschenburg examples
themselves. Let $U_{klpq}$ be the one-parameter subgroup of $SU(3)\times
SU(3)$  from   Example 4.6 and endow  $SU(3)$  with  a bi-invariant
metric. We will show below that for all $m$,  the geodesic flow of the
Eschenburg manifold  $E_m \df SU(3) / U_{1,-1,2m,2m}$ endowed with the
submersion metric is completely integrable.  As Eschenburg showed,  this
is another example of  a strongly inhomogeneous manifold
\cite{Eschenburg}. Also notice that for $m=0$ we obtain a Wallach manifold
\cite{Wallach}. 

\def \U{U_{1,-1,2m,2m}}
\def \r{{\cal R}}
\def \v{{\cal H}}
\def\tr{{\mbox tr }}
For simplicity, set $U=\U$.
 Denote by $\Phi :T SU(3) \rightarrow su(3)+su(3)$ the moment map of the
action of $SU(3)\times SU(3)$ on the tangent bundle $TSU(3)$. 
As in Example 4.6, we  use the   ascending chain of subgroups
$\{1\}\times U(1)\subset U(1)\times U(1)\subset U(1)\times  U(2)\subset
U(2)\times U(2)\subset U(2)\times SU(3)\subset SU(3)\times SU(3)$. Let
$pr_1$ and $pr_2$ be the projections of $su(3) +su(3)$  onto the  first
and   second factor respectively. Denote by $pr_{u(i)}$ the  orthogonal
projection of $su(3)$ to $u(i)$.  Further identify $su(3) ^*$ with $su(3)$
via the Cartan-Killing form as usual.  Then the Thimm functions on
$TSU(3)$ are the pull backs under the moment map of the following
functions on $su(3) +su(3)$:
\[ \begin{array}{ll}
f_1 = i\:\tr (\xi) \circ pr_{u(1)} \circ  pr_1\:\:\:\: \:\:\: \:\: &
\:\:\:\:\:\: f_5 = i\:\tr (\xi ^3 ) \circ pr_1 \\
f_2 = i\:\tr (\xi) \circ pr_{u(2)} \circ pr_1 \:\:\:\:\:\:\: \:\: &\:\:
\:\:\:\: f_6 = i\: \tr (\xi) \circ pr_{u(1)} \circ  pr_2 \\
 f_3 =\tr (\xi ^2)  \circ pr_{u(2)} \circ pr_1 \:\:\:\: \:\: \:\:\:&
\:\:\:\:\:\:   f_7 =i \:\tr (\xi) \circ pr_{u(2)} \circ pr_2\\
 f_4 = \tr (\xi ^2)  \circ pr_1 \:\:\:\: \:\:\: \:\: & \:\:\:\:\:\:  f_8
=\tr (\xi ^2 ) \circ pr_{u(2)} \circ pr_2.
\end{array}\]
As in   Example 4.6, all the Thimm integrals  $f_i \circ \Phi$ on $TSU(3)$
are invariant under $U$, and thus induce Poisson-commuting functions
$\tilde{f_i}$  on  $TE_m$.

Let us now show the  independence of  seven of these functions, namely
$\tilde{f_2}, \ldots, \tilde{f_8}$. First note that by  real analyticity
we only need to establish  the independence  of these functions at one
point.

% For simplicity, set $U = \U$. Denote by  $\Phi :T SU(3) \rightarrow
% su(3) + su%(3)$  the moment map of the action of $SU(3) \times SU(3)$ on
% the  tangent bund%le $T SU(3)$. 
      Let $\v$ denote the set of horizontal vectors on $TSU(3)$. Recall
that $\v = \Phi _U ^{-1} (0)$ where $\Phi _U$ is the moment  map of the
action of $U$ on $TSU(3)$.

First  we will reduce the problem to a calculation in the Lie algebra.
Suppose that the  $\tilde{f_i} $, $i=2,\ldots,8$, are dependent at the
projection  $v$ of  a vector $\hat{v} \in \v$ via some relation $\sum
_{i=2} ^8  c_i d\tilde {f_i} =0$ on $T_v E_m$. Then the 1-form $\sum
_{i=2} ^8 c_i d(f_i \circ \Phi)$ is 0 on horizontal lifts of double
tangent vectors. 
 Since the functions $ f_i \circ \Phi$ are $U$-invariant, $\sum _{i=2} ^8
c_i d(f_i \circ \Phi)$ is also 0 on tangent vectors to the $U$-orbit of
$\hat{v}$. This implies that $\sum _{i=2} ^8 c_i d(f_i \circ \Phi) =0$ on
$T_{\hat{v}} \v$.  Now suppose that $\Phi (\v)$ is a manifold in a
neighborhood of $\Phi (\hat{v})$ and that  $\Phi (\hat{v})$ is a regular
value of $\Phi:\v \rightarrow \Phi (\v)$.
 Then a dependence of the restrictions of the $f_i \circ \Phi$ to $\v$ at
$\hat{v}$ implies a dependence of the restrictions of the $f_i$ to the
image of $\v$ under $\Phi$ at $\Phi (\hat{v})$. This is the reduction to a
calculation in the Lie algebra.

Next we need to determine $\Phi (\v)$. First let us  describe $\Phi$
itself. 
 As usual identify $T_1 SU(3)$ with the orthogonal complement (with
respect to the Cartan-Killing form) of the diagonal embedding $\Delta
su(3)$ of $su(3)$ into $su(3) \times su(3)$, that is with 
$\{(X,-X) \mid X \in su(3)\}$. Then we have the following formula for the
value of the moment map at a translate of a vector  $(X,-X)$ in $T_1
SU(3)$
\[ \Phi ((g_1,g_2) _* (X,-X)) = (\Ad g_1 (X), -\Ad g_2 (X)).\]
Thus $\r \df \Phi (TSU(3))$ is given by
\[ \r = \{(X,Y) \mid X \mbox{ is conjugate to } -Y \mbox{ in }su(3)\}.\]
Next note that $\Phi _U = i^* \circ \Phi$ where $i: u \rightarrow su(3)
+su(3)$ is the embedding of the Lie algebra $u$ of $U$ into $su(3)
+su(3)$. Thus
 the horizontal vectors in $TSU(3)$ are the preimage $\v = \Phi _{U} ^{-1}
(0)=\Phi ^{-1} (u^{\perp})$. Hence the image of $\v$ under $\Phi$ is
\[ \Phi (\v) = \r \cap u^{\perp} .\]

Let us now give an outline of the calculations that show  the independence
of the restrictions of  $f_2, \ldots, f_8$ to $\Phi(\v)$ at the   point $p
\in \Phi (\v)$ given by $p = (P,-P)$ where 
\[ P = \left(\begin{array}{ccc}
0  &  2  &  1  \\
-2 &  0  &  0  \\
-1 &  0  &  0
\end{array} \right).\]
 One easily shows that $ \r \cap u^{\perp}$ is a manifold in a
neighborhood of $p$, and that $p$ is a regular value of $\Phi$. Suppose
that on $T_p (\r \cap u^{\perp})$ we have
\[\:\:\: (*) \:\:\: \:\:\: \:\:\:  \:\:\:  \:\:\:  \:\:\: \:\:\:  \:\:\:
\:\:\: \:\:\: \:\:\: \:\:\: \:\:\:  \sum _{i=2} ^8  c_i d {f_i} =0 .
\:\:\:  \:\:\: \:\:\: \:\:\: \:\:\: \:\:\: \:\:\: \:\:\: \:\:\: \:\:\:
\:\:\: \:\:\: \:\:\: \]
We will exhibit several tangent vectors in $T_p (\r \cap u^{\perp})$ which
force various relations between the coefficients $c_i$, forcing them to be
0 eventually. 

\begin{enumerate}
\item Set $p _1 ^t = (P _1 ^t,-P _1 ^t)$ where 
\[ P _1 ^t = \left(\begin{array}{ccc}
0  &  2  &  1 +t \\
-2 &  0  &  0  \\
-1-t &  0  &  0
\end{array} \right).\]
Then $p _1 ^t \in  \r \cap u^{\perp}$. Since the nontrivial projections
$pr _{u(i)}$ of $P _1 ^t$ are all constant, only $d f_4$ and $df_5$ can be
nonzero on $ v_1=\frac{d}{dt} \Big| _{t=0} \, p_1 ^t$. The eigenvalues of
$P _1 ^t$ are $0,\sqrt{-4-(1+t)^2}$ and $-\sqrt{-4-(1+t)^2}$. Therefore
we get $df_5 (v_1)=0$, $df_ 4 (v_1) \neq 0$, and thus  $c_4 =0$.
\item Set $p _2 ^t = (P _2 ^t,-P _2 ^t)$ where 
\[ P _2 ^t = \left(\begin{array}{ccc}
0  &  2  &  1  \\
-2 &  0  &  it  \\
-1 &  it  &  0
\end{array} \right).\]
Then $p _2 ^t \in  \r \cap u^{\perp}$. The eigenvalues of $P _2 ^t$
satisfy the equation
\[- \lambda ^3 -\lambda (t^2 +5) -4it =0.\]
Hence $f_5 (p _2 ^t) = - i(\tr (P _2 ^t) (t^2 +5) -12 it)=-12 t$.
 As above and since $c_4 =0$ we conclude that $c_5 =0$.
\item  Let $s(t) = \sqrt{(t+2)^2 -3}$ and set $p _3 ^t = (P _3 ^t,-Q _3
^t)$ where 
\[ P _3 ^t = \left(\begin{array}{ccc}
0  &  2 +t &  1  \\
-2 -t &  0  & 0  \\
-1 &  0  &  0
\end{array} \right)
\mbox{  and } 
Q _3 ^t = \left(\begin{array}{ccc}
0  &  2  &  s(t)\\
-2  &  0  & 0  \\
-s(t) &  0  &  0
\end{array} \right).\]
A calculation of the eigenvalues shows that $p _3 ^t \in \r \cap
u^{\perp}$ and that only $df_3$ gives a nonzero contribution in  (*) when
applied to  $\frac{d}{dt} \Big| _{t=0} \, p_3 ^t$. Therefore we get
$c_3=0$.
\item Considering $p _4 ^t = (Q _3 ^t,-P _3 ^t)$ we find that $c_8 =0$.
\item Set  $p_5 ^t = (P,-\Ad (\exp t A)(P))$ where 
\[ A = \left(\begin{array}{ccc}
0  & i  & 0 \\
i  & 0  & 0 \\
0  & 0  & 0 
\end{array}\right) .\]         
Clearly, $p_5^t$ lies in ${\cal R} \cap u ^{\perp}$ since $id \times 
A$ commutes with $u$. Since 
\[[A,P] = \left(\begin{array}{ccc}
-4i  &  0  &  0 \\
0    & 4i  &  i  \\
0    &  i  &  0
\end{array}\right), \]                                            
 we get $df_2 (\frac{d}{dt} \Big| _{t=0} p_5 ^t) = df_7 (\frac{d}{dt} 
\Big| _{t=0} p_5 ^t)  =0$ while $df_6 (\frac{d}{dt} \Big| _{t=0} p_5 ^t)
\neq 
0$. This implies $c_6 =0$.
\item Let $p_6 ^t = (P_6 ^t, -P_6 ^t)$ where 
\[ p_6 ^t = \left(\begin{array}{ccc}
t  &  2                     &  1  \\
-2 &  -\frac{6m-1}{6m+1} t  &  0  \\
-1 &  0  & \frac{-2}{6m+1} t
\end{array}\right).\]
One sees easily that $p_6^t \in {\cal R} \cap u ^{\perp}$, and then 
that $c_2 = c_7$.
\item Suppose that $c_2 = c_7 \neq 0$. Then $df_2 = -df_7$ on $T_p ({\cal
R} \cap u ^{\perp})$. Note that 
\[ \left( \exp t \left(\begin{array}{ccc}
0  &  i  &  0  \\
i  &  0  &  0  \\
0  &  0  &  0
\end{array}\right) P, -P  \right)\in {\cal R}.\]
The tangent vector $v$ to this  curve at 0 is given by 
\[\left( \left[ \left(\begin{array}{ccc}
0  &  i  &  0  \\
i  &  0  &  0  \\
0  &  0  &  0
\end{array}\right) , P \right] , 0\right) =
\left( \left(\begin{array}{ccc}
-4i &  0  &  0  \\
0   &  4i &  i \\
0   &  i  &  0
\end{array}\right) ,0 \right) . \]
Thus $v$ is not perpendicular to $u$ while $df_2 (v) =df_7(v)=0$.  Hence
$df_2 =- df_7$ on $T_p ({\cal R} \cap u ^{\perp}) + \R v = T_p {\cal R}$. 

On the other hand, consider the curve in ${\cal R}$ given by 
\[ \left( \exp t \left(\begin{array}{ccc}
0  &  0  &  i  \\
0  &  0  &  0  \\
i  &  0  &  0
\end{array}\right) P, -P  \right).\]
Its tangent vector  $w$ at 0 is 
\[ \left( \left[ \left(\begin{array}{ccc}
0  &  0  &  i  \\
0  &  0  &  0  \\
i  &  0  &  0
\end{array}\right) , P \right] , 0\right) =\left( \left(\begin{array}{ccc}
-2i &  0  &  0  \\
0   &  0  &  2i \\
0   &  2i &  2i
\end{array}\right) ,0 \right) . \]
Clearly we have  $df_2 (w) \neq df_7 (w)$, a contradiction. Therefore we
get $c_2 = c_7 =0$, and $f_2, \ldots, f_8$ are a.e. independent.
\end{enumerate}
}\end{Example}

As a final application of the submersion method we construct a Riemannian
metric with completely integrable geodesic flow on an exotic sphere. Again
the submersion in question does not have the product type. The integrals
themselves arise both from the submersion method combined with a Thimm
construction as well as from the isometry group of this exotic sphere.  

\noindent {\bf Example 4.8} Consider the exotic 7-sphere $\Sigma$
constructed by Gromoll and Meyer in \cite{Gromoll-Meyer}. It arises as a
biquotient of $Sp(2)$ by the following action of $Sp(1)$. For $q \in
Sp(1)$ and  $Q \in Sp(2)$ set 
\[ (q,Q) \mapsto 
\left( \begin{array}{cc}
q & 0 \\
0 & 1
\end{array} \right)
Q
\left( \begin{array}{cc}
\bar{q} & 0 \\
0 & \bar{q}
\end{array} \right) \]
where $\bar{q}$ denotes the complex conjugate of $q$. This also defines an
embedding $U$ of $Sp(1)$ into $Sp(2) \times Sp(2)$. Note that the
canonical $O(2)$ in $Sp(2)$ commutes with the right action of $Sp(1)$
while an obvious $Sp(1)$ commutes with the left action. We give $\Sigma$
the submersion metric determined by the biinvariant metric on $Sp(2)$. 

The basic argument is much the same as in Example 4.7. Again, let $\Phi:
TSp(2) \rightarrow sp(2) + sp(2)$ denote the moment map of the action of
$Sp(2) \times Sp(2)$ on the tangent bundle $TSp(2)$. Let $pr_1$ and $pr_2$
denote the projections of $sp(2) + sp(2)$ onto the first and second factor
respectively. Further we denote the orthogonal projection to a subalgebra
$h \subset sp(2)$ by $pr_{h}$. We embed $sp(2)$ into $u(4)$ canonically.
Then we define  the following functions on $sp(2) + sp(2)$ using complex
valued traces:
\[\begin{array}{ll}
f_1= \tr (\xi ^2) \circ pr_{sp(1) \times 1} \circ pr_1
	& f_5=\tr (\xi ^4) \circ pr_1 \\
f_2= \tr (\xi ^2) \circ pr_{so(2)} \circ pr_2                         
	& f_6= \tr (\xi ^2) \circ pr_{1 \times l} \circ pr_1      \\
f_3=	\tr (\xi ^4) \circ pr_{sp(1) \times sp(1)} \circ pr_1
& f_7=\tr (\xi ^2) \circ pr_{sp(1) \times sp(1)} \circ pr_2 \\
f_4=\tr (\xi ^2) \circ pr_1
\end{array}\]
where $so(2)$ refers to the Lie algebra of the canonical $O(2)$  above
while $l$ refers to  the subalgebra of $sp(1)$ generated by 
\[ \left(\begin{array}{cc}
0   &  0 \\
0   &  k
\end{array} \right).\]
 These functions are all invariant functions on  some subalgebra pulled
back to $sp(2) + sp(2)$. These subalgebras are either contained in each
other as in Thimm's argument or commute with each other. It easily follows
that they all Poisson commute. Note that $f_2$ and $f_6$ are just first
integrals coming from the the isometry group. 
Also note that all these functions are invariant under the adjoint action
of $U$. Hence their pullbacks to $TSp(2)$ under $\Phi$ are invariant under
the action of $Sp(1)$ on $Sp(2)$, and thus they descend to functions
$\tilde{f}_i$, $i=1,\ldots,7$ on $T \Sigma$. As in Example 4.7 the
independence of the $\tilde{f}_i$ at the projection of a horizontal vector
$\hat{v}$  is equivalent to the independence of the restrictions of
$f_1,\ldots, f_7$ to ${\cal R} \cap u^{\perp}$ near $\Phi (\hat{v})$ where
$u$ is the  Lie algebra of $U$ and ${\cal R} = \{(X,Y) \mid X \mbox { is
conjugate to } -Y\}$. We assume here that $\Phi (\hat{v})$ is a regular
value of $\Phi : {\cal H} \rightarrow \Phi ({\cal H})$ and that ${\cal R}
\cap u^{\perp}$ is a manifold near $\Phi (\hat{v})$. 

Next we will indicate a point $p$ in ${\cal R} \cap u^{\perp}$ and tangent
vectors in $T_p ({\cal R} \cap u^{\perp})$ that show the independence of
$f_1,\ldots, f_7$. 

Let 
\[ F(t) \df \left( \begin{array}{cc}
\cos \pi t         &   \sin \pi t   \\
- \sin \pi t       &   \cos \pi t
\end{array}\right) \]
and 
\[ G(t) \df \left( \begin{array}{cc}
1          &             0   \\
0          &  \cos \pi t + j \sin \pi t 
\end{array} \right)\]
where $1,i,j$ and $k = ij$ are the standard basis of the quaternions. 
Set $Q \df F(\frac{1}{3}) G(\frac{1}{4})$ and 
\[ P \df \left( \begin{array}{cc}
2i -2j - \frac{149 +18 \sqrt{2} \sqrt{3}}{9} k     &
              1 + 3i +2j -3k           \\
-1 +3i+2j-3k        &
		  5i+(6+\sqrt{2} \sqrt{3})j +\frac{26}{3} k
\end{array}\right). \]
Define $R = Q P Q^*$ where $Q^*$ is the conjugate transpose of $Q$, and
set $p =(R,P)$. Then $p \in {\cal R} \cap u^{\perp}$. One can check that
$p$ is a regular value of $\Phi$, that ${\cal R} $ and $u^{\perp}$
intersect transversally at $p$ and that ${\cal R} \cap u^{\perp}$ is a
manifold near $p$. 

Next we will list the relevant tangent vectors in $T_p ({\cal R} \cap
u^{\perp})$. We need the following matrices:
\[ \begin{array}{cccc}
D_1 \df  \left( \begin{array}{cc}
    i   &  0 \\
    0   &  0
\end{array} \right) \:\:\:\: &\:\:\:\:
D_2 \df \left( \begin{array}{cc}
    0   &  0 \\
    0   &  i
\end{array} \right) \:\:\:\: &\:\:\:\:
D_3 \df \left( \begin{array}{cc}
    0   &  j \\
    j   &  0
\end{array} \right) \:\:\:\: &\:\:\:\:
D_4 \df \left( \begin{array}{cc}
    0   &  k \\
    k   &  0
\end{array} \right)      
\end{array} \]  
\[ \begin{array}{cccc}
D_5 \df \left( \begin{array}{cc}
    0   &  1 \\
    -1   &  0
\end{array} \right)  \:\:\:\:\: &\:\:\:\:\:
D_6 \df \left( \begin{array}{cc}
    k   &  0 \\
    0   &  0
\end{array} \right) \:\:\:\:\: &\:\:\:\:\:
D_7 \df \left( \begin{array}{cc}
    0   &  0 \\
    0   &  j
\end{array} \right) . 
\end{array} \]
Then the tangent vectors are:
\begin{enumerate}
\item $v_1 =(D_5,- Q^* D_5 Q)$
\item $v_2 = (R,P)$
\item $v_3 = (0,[D_5,P])$
\item $v_4 = (- Q (\sqrt{2} \sqrt{3} D_7 + D_3) Q^*,\sqrt{2} \sqrt{3} D_7
+ D_3)$
 \item $v_5=(-Q D_1 Q^*,D_1 +[\frac{3}{80} D_4 -\frac{9}{80} D_6,P])$
\item $v_6=([D_7,R],0)$
\item $v_7 = (-D_2,Q^* D2 Q +[\frac{1}{48} D_3 +\frac{1}{40} D_4
-\frac{3}{40} D_6,P])$.
\end{enumerate}
As in Example 4.7, evaluating  a relation of the $df_i$'s on $T_p ({\cal
R} \cap u^{\perp})$ on these seven tangent vectors forces this relation to
be trivial. This long calculation as well as finding the vectors above was
done  by computer using Mathematica.

\section{\em Connected Sums.}   
{\rm In this section we will combine the submersion technique of the last
section with a  glueing trick to construct metrics on $\C P^{n+1}\#\C
P^{n+1}$ for $n$ even with completely integrable geodesic flows.
 Topologically these spaces are obtained from two copies of
$S^{2n+1}\times_{S^{1}}D^{2}$ where $D^{2}$ is the 2-disk and $S^{1}$ acts
diagonally, glued along their boundary
$S^{2n+1}\times_{S^{1}}S^{1}=S^{2n+1}$ by an orientation reversing map.
 The metrics that we will use were already considered in \cite{C}. Let us
describe them.

Consider the Hopf fibration $S^{1}\rightarrow S^{2n+1}\rightarrow \C
P^{n}$ and endow $S^{2n+1}$ with the metric $g_{t}$ as in Example 4.4.
 Now equip $\R^{2}$ with a metric $h_{t}$ ($t^{2}\neq 1$) given in polar
coordinates by :
\[h_{t}(\partial/\partial r,\partial/\partial
r)=1\;\;\;\;h_{t}(\partial/\partial r,\partial/\partial
\theta)=0\;\;\;\;
h_{t}(\partial/\partial\theta,\partial/\partial\theta)=f^{2}_{t}(r)\]
where $f_{t}(r)$ is a smooth function with the properties $f_{t}(0)=1$,
$f_{t}'(0)=1$ and $f_{t}(r)\equiv 2\pi t^{2}/\sqrt{t^{4}-1}$ for
sufficently big $r>R$.  

Set $\eta=S^{2n+1}\times_{S^{1}}\R^{2}$ with the submersion metric.
 If we restrict to the disk bundle $D_{\bar{R}}(\eta)$ with $\bar{R}>R$,
then an annular neighborhood of the boundary splits isometrically as
$\partial D_{\bar{R}}(\eta)\times I$ where $I$ denotes an interval.
 In fact, $A=\{X \in \R ^{2} \; \mid R<\parallel X\parallel< \bar{R} \}$
splits isometrically as $S^{1}\times I$ and $S^{1}$ acts trivially on $I$.
 Then 
\[S^{2n+1}\times_{S^{1}}A=S^{2n+1}\times_{S^{1}}(S^{1}\times
I)=(S^{2n+1}\times_{S^{1}}S^{1})\times I=S^{2n+1}\times I\]
and $S^{2n+1}=\partial D_{\bar{R}}(\eta)$ gets back the metric of constant
curvature.
Since the metric splits as a product $S^{2n+1}\times I$ near the boundary,
by glueing two such disk bundles we get a smooth metric on $\C P^{n+1} \#
\C P^{n+1}$.

According to Example 4.4 the metric on the disk bundle $D_{\bar{R}}(\eta)$
is completely integrable with first integrals
$f_{1},...,f_{2n+1},f_{2n+2}$.
 In fact $f_{1},...,f_{2n+1}$ are induced by the Thimm integrals on the
tangent bundle of $(S^{2n+1},g_{t})$ and $f_{2n+2}$ is induced by the
metric $h_{t}$ (cf. Proposition 4.2).
 All the $f_{i}\;'s$ are invariant under derivatives of translations on
$I$. Therefore they will fit together smoothly with the integrals on the
second $D_{\bar{R}}(\eta)$ if they happen to be invariant under the
derivative of the orientation reversing map that we use for the glueing.
 
As a  glueing map on the boundary $S^{2n+1}$ we will take the complex
conjugation $\tau$ i.e. the restriction to $S^{2n+1}\subset \C^{n+1}$ of
the map:
\[(z_{1},...,z_{n+1})\rightarrow (\bar{z}_{1},...,\bar{z}_{n+1})\]
This map is orientation reversing for $n$ even (for $n$ odd, one
rediscovers Example 4.4).

As we will see below some of the functions $f_{i}$ are not invariant under
$d\tau$. Thus a small modification will be needed.

Denote by $\pi$ the projection map $\pi: S^{2n+1}\times S^{1}\rightarrow
S^{2n+1}\times_{S^{1}}S^{1}$ and by $\sigma$ the map $\sigma
:S^{2n+1}\times S^{1}\rightarrow S^{2n+1}\times S^{1}$ given by
$\sigma(z_{1},...,z_{n+1},e^{i\theta})=
(\bar{z}_{1},...,\bar{z}_{n+1},e^{-i\theta})$.
 Note that $\sigma$ takes $S^{1}$-orbits into $S^{1}$-orbits since
\[e^{i\varphi}.\sigma(z_{1},...,z_{n+1},e^{i\theta})=
\sigma(e^{-i\varphi}.(z_{1},...,z_{n+1},e^{i\theta}))\]
Hence $\sigma$ descends to a map
$\hat{\sigma}:S^{2n+1}\times_{S^{1}}S^{1}\rightarrow
S^{2n+1}\times_{S^{1}}S^{1}$.
 Observe that under the natural diffeomorphism $\psi:S^{2n+1}\rightarrow
S^{2n+1}\times_{S^{1}}S^{1}$ given by
$\psi(z_{1},...,z_{n+1})=\pi(z_{1},...,z_{n+1},1)$, the map $\hat{\sigma}$
is complex conjugation, i.e. $\psi^{-1}\circ\hat{\sigma}\circ \psi=\tau$.
% In fact let us compute:
%\[\psi^{-1}\circ
%\hat{\sigma}\circ\psi(z_{1},...,z_{n+1})=
%\psi^{-1}\circ\hat{\s%igma}\circ\pi(z_{1},...,z_{n+1},1)=\]
%

%\[=\psi^{-1}\circ\pi\circ\sigma(z_{1},...,z_{n+1},1)=
%(\bar{z}_{1},...,\bar{z}_{%n+1})\]

Since $h_{t}$ is invariant under the map $(r,\theta)\rightarrow
(r,-\theta)$ we deduce that the integral $f_{2n+2}$ will be invariant
under the derivative of $\tau$.

Therefore we need to find integrals on $(S^{2n+1},g_{t})$ which are
invariant under $d\tau$ and under the $S^{1}$-action. In view of the
previous arguments this automatically implies that the induced integrals
on $S^{2n+1}\times_{S^{1}}S^{1}$ are also invariant under $d\tau$ and that
we will be able to fit them smoothly.

Recall that the integrals we have on $(S^{2n+1},g_{t})$ were obtained by
the Thimm method using the action of the group $SU(n+1)\times S^{1}$. 
 Let $f_{2n+1}$ denote the integral induced by the $S^{1}$-action.
 Since
$e^{i\varphi}.\tau(z_{1},...,z_{n+1})=\tau(e^{-i\varphi}.(z_{1},...,z_{n+1}))$,
we see that $f_{2n+1}$ is not invariant under $d\tau$. But $f^{2}_{2n+1}$
is clearly invariant and still is a first integral.
 We will now use a similar trick for the integrals that arise from the
$SU(n+1)$-action.

Identify $S^{2n+1}$ with $SU(n+1)/SU(n)$ in the usual way, i.e. by means
of the diffeomorphism $[A]\rightarrow  A(1,0,...,0)$ where $[A]$ denotes
the equivalence class of a matrix $A\in SU(n+1)$. 
 Since $\tau\circ A=\bar{A}\circ\tau$ it is easy to check that $\tau$
operates on $SU(n+1)/SU(n)$ as the map $[A]\rightarrow [\bar{A}]$.
 Decompose $su(n+1)$ as $su(n)\oplus m$ where $m$ denotes the orthogonal
complement of $su(n)$ in $su(n+1)$ with respect to the standard Killing
form. The moment map $\phi$ of the action of $SU(n+1)$ on the tangent
bundle of $SU(n+1)/SU(n)$ can be written as (cf. \cite[Lemma 3.2]{Thimm}):
\[\phi(dL_{A}(B))=\Ad_{A}(B)\]
where  $A\in SU(n+1)$, $B\in m$ and $L_{A}$ denotes the left translation
on $SU(n+1)/SU(n)$.

Since the integrals arising from the Thimm method have the form
$h\circ\phi$ where $h\in C^{\infty}(su(n+1))$, they are invariant under
the derivative of $\tau[A]=[\bar{A}]$ if and only if for every $A\in
SU(n+1)$ and $B\in m$ we have
\[h(\Ad_{\bar{A}}(\bar{B}))=h(\Ad_{A}(B)).\]
%But $\Ad_{\bar{A}}(\bar{B})=\bar{\Ad _{A}(B)}$.
 Therefore $h\circ\phi$ is invariant under $d\tau$ if $h$ is invariant
under conjugation on $su(n+1)$.
 If $B\in su(n+1)$ then $\bar{B}=-B^{t}$. Hence  we need $h$ such that
$h(B)=h(-B^{t})$.

Denote by $\pi_{j}:su(n+1)\rightarrow u(j)$ the map defined by:
\[su(n+1)\ni
\left(\begin{array}{cl}\alpha&\beta\\\gamma&\delta\end{array}\right)\rightarrow
\delta\in u(j)\]
The 2n-functions in involution on $su(n+1)$ that we get from the Thimm
method are (\cite[Proof of Theorem 7.4]{Thimm}):
\[h_{j}(B)=-\frac{1}{2}itr(\pi_{j}B)\;\;\;\;\;j=1,...,n\]
\[h_{n+j-1}(B)=-\frac{1}{4}tr(\pi_{j}B)^{2}\;\;\;\;\;j=2,...,n+1\]
Clearly the $h_{n+j-1}\;'s$ are invariant under $B\rightarrow -B^{t}$, but
the $h_{j}\;'s$ are not. Instead consider the functions:
\[h_{j}^{2}(B)\;\;\;\;\;j=1,...,n\]
\[h_{n+j-1}(B)\;\;\;\;\;j=2,...,n+1\]
Now they are all invariant under $B\rightarrow -B^{t}$, they are still in
involution and they are functionally independent a.e.
 Hence the pull back of these functions by the moment map $\phi$ gives a
set of 2n-functions that verifies all the  necessary conditions.}

\end{document}